\documentclass{article} 

\usepackage{amssymb}
\usepackage{eucal}
\usepackage{amsmath}

\usepackage{graphicx}
\usepackage{float}
\usepackage{framed}
\usepackage{blkarray}
\usepackage{enumitem}

\usepackage{relsize}

\newtheorem{thm}{\bf Theorem}[section]
 
\newtheorem{prop}[thm]{\bf Proposition}

\usepackage{graphicx}

\newcommand{\wt}{\widetilde}

\newcommand{\R}{{\mathbb  R}} 
\DeclareMathOperator{\tr}{tr}

\begin{document}

\title{Second order optimality on orthogonal Stiefel manifolds}

\author{Petre Birtea, Ioan Ca\c su, Dan Com\u{a}nescu\\
{\small Department of Mathematics, West University of Timi\c soara} 
\\
{\small Bd. V. P\^ arvan, No 4, 300223 Timi\c soara, Rom\^ania}\\
{\small Email: petre.birtea@e-uvt.ro, ioan.casu@e-uvt.ro, dan.comanescu@e-uvt.ro}}
\date{}

\maketitle

\begin{abstract}
The main tool to study a second order optimality problem is the Hessian operator associated to the cost function that defines the optimization problem. 
By regarding an orthogonal Stiefel manifold as a constraint manifold embedded in an Euclidean space we obtain a concise matrix formula for the Hessian of a cost function defined on such a manifold. We introduce an explicit local frame on an orthogonal Stiefel manifold in order to compute the components of the Hessian matrix of a cost function. We present some important properties  of this frame. 
As applications we rediscover second order conditions of optimality for the Procrustes and the Penrose regression problems (previously found in the literature). For the Brockett problem we find necessary and sufficient conditions for a critical point to be a local minimum. Since many optimization problems are approached using numerical algorithms, we give an explicit description of the Newton algorithm on orthogonal Stiefel manifolds.
\end{abstract}

\section{Introduction}

Optimization problems on Stiefel manifolds appear in important applications such as statistical analysis of data \cite{turaga}, blind signal separation \cite{lu}, distance metric learning \cite{shukla}, among many other problems.

In Section \ref{section-2} we regard an orthogonal Stiefel manifold as a preimage of a regular value for a set of constraint functions.
We adapt the method presented in the papers \cite{birtea-comanescu} and \cite{Birtea-Comanescu-Hessian} (the so called embedded gradient vector field method) to the particular case of Stiefel manifolds. We embed the Stiefel manifold $St^n_p$ in the larger Euclidean space $\mathcal{M}_{n\times p}(\R)$, in order to take advantage of the simpler geometry of this Euclidean space. This setting allows us to present necessary and sufficient conditions for critical points of a cost function defined on a Stiefel manifold and a formula for the Hessian of this cost function in a concise matrix form. This formula is important in the study of the second order optimality.

In order to explicitly compute the components of the Hessian matrix of a cost function, in Section \ref{section-3} we introduce an explicit local frame for an orthogonal Stiefel manifold and we present some important properties of this local frame. We determine the components of the Hessian matrices of the constraint functions in this frame.

In the last section we apply the results of the previous sections to some important problems arising from practical applications. For the Procrustes and the Penrose regression problems we obtain second order conditions, which have been previously presented in \cite{chu} using a different approach, which involves the projected Hessian. For the Brockett problem, see \cite{absil-mahony-sepulchre-1}, we find necessary and sufficient conditions for a critical point of the cost function to be a local minimum. As an example, for a particular Brockett cost function defined on the orthogonal Stiefel manifold  $St^4_2$ we give a  list of the critical points and we completely characterize them. 

In many cases optimization problems are approached using numerical methods on manifolds. The Newton algorithm is a type of algorithm that uses the second order information about a cost function. We give an explicit description of the Newton algorithm in the case of an orthogonal Stiefel manifold. There exists a rich literature that deals with the construction of Newton algorithm on manifolds, see \cite{absil-mahony-sepulchre-1}, \cite{gabay}, \cite{manton}, and \cite{wen}.

\section{Hessian matrix on orthogonal Stiefel manifolds}\label{section-2}

Let $\mathcal{S}\subset \mathfrak{M}$ be a submanifold of a Riemannian manifold $(\mathfrak{M},{\bf g})$, which can be described by a set of constraint functions, i.e. $\mathcal{S}={\bf F}^{-1}(c_0)$, where ${\bf F}=(F_1,\dots,F_k):\mathfrak{M}\rightarrow \R^k$ is a smooth map and $c_0\in \R^k$ is a regular value of ${\bf F}$. The manifold $\mathcal{S}$ becomes a Riemannian manifold when endowed with the induced metric ${\bf g}_{_{ind}}$. 

The Riemannian geometry of the submanifold can be more complicated than the Riemannian geometry of the ambient manifold.
In optimization problems we need, in general, to compute the gradient vector field and the Hessian operator of a smooth cost function $\widetilde{G}:(\mathcal{S},{\bf g}_{_{ind}})\rightarrow \R$. In what follows  we present a method, called  {\bf the embedded gradient vector field method},  to compute the gradient vector field and the Hessian operator of $\widetilde{G}$ using only the geometry of the ambient manifold $(\mathfrak{M},{\bf g})$.

Let $G:(\mathfrak{M},{\bf g})\rightarrow \R$ be a smooth prolongation of $\widetilde{G}$. In \cite{birtea-comanescu}, \cite{Birtea-Comanescu-Hessian}, \cite{5-electron}, it has been proved that
\begin{equation}\label{partial-G-134}
\nabla _{{\bf g}_{_{ind}}}{G}_{|{\bf F}^{-1}(c)}(s)=\partial_{\bf g} G(s),\,\,\,\forall s\in {\bf F}^{-1}(c),\,\,\,\text{and}\,\,c\,\,\text{an arbitrary regular value},
\end{equation}
where $\partial_{\bf g} G$ is defined on the open set of regular points $\mathfrak{M}^{reg}\subset \mathfrak{M}$ of the constraint function, and it is the unique vector field that is tangent to the foliation generated by ${\bf F}$ having property \eqref{partial-G-134}. {The embedded gradient vector field}  is given by the following formula:
\begin{equation*}
\partial_{\bf g} G(s)=\nabla_{\bf g} G(s)-\sum\limits_{i=1}^k\sigma_{\bf g}^{i}(s)\nabla_{\bf g}F_i(s),\,\,\,\forall\,s\in \mathfrak{M}^{reg}.
\end{equation*}
The Lagrange multiplier functions $\sigma_{\bf g}^{i}:\mathfrak{M}^{reg}\rightarrow \R$ are defined by the formula

\begin{equation*}\label{sigma-101}
\sigma^i_{\bf g}(s):=\frac{\det \left(\text{Gram}_{(F_1,\ldots ,F_{i-1},G, F_{i+1},\dots,F_k)}^{(F_1,\ldots , F_{i-1},F_i, F_{i+1},...,F_k)}(s)\right)}{\det\left(\text{Gram}_{(F_1,\ldots ,F_k)}^{(F_1,\ldots ,F_k)}(s)\right)},
\end{equation*}
where
\begin{equation*}\label{sigma}
\text{Gram}_{(g_1,...,g_s)}^{(f_1,...,f_r)}:=\left[%
\begin{array}{cccc}
  {\bf g}(\nabla_{\bf g} g_1,\nabla_{\bf g}f_{1}) & ... & {\bf g}(\nabla_{\bf g} g_s,\nabla_{\bf g} f_{1}) \\
  \vdots & \ddots & \vdots \\
  {\bf g}(\nabla_{\bf g} g_1,\nabla_{\bf g}f_r) & ... & {\bf g}(\nabla_{\bf g} g_s,\nabla_{\bf g} f_r) 
\end{array}%
\right].
\end{equation*}
As mentioned in \cite{Birtea-Comanescu-Hessian} in a critical point $s_0\in \mathcal{S}$ of $\widetilde{G}$  the numbers $\sigma^i_{\bf g}(s_0)$ coincide with the classical Lagrange multipliers.

Also, in \cite{Birtea-Comanescu-Hessian}, \cite{5-electron} it has been proved that 
\begin{equation*}\label{Hessian-general}
\text{Hess}_{{\bf g}_{_{ind}}}\, \widetilde{G}(s) =\left(\text{Hess}_{\bf g}\,G(s)-\sum_{i=1}^k\sigma_{\bf g}^{i}(s)\text{Hess}_{\bf g}\, F_i(s)\right)_{|T_s \mathcal{S}\times T_s \mathcal{S}},\,\,\,\forall s\in \mathcal{S}.
\end{equation*}

We apply the above general setting to the particular case of the orthogonal Stiefel manifold $St_p^n$, regarded as a submanifold  of $\mathcal{M}_{n\times p}(\R)$.
For $n\geq p\geq 1$, we consider the orthogonal Stiefel manifold:
$$St_p^n=\{U\in \mathcal{M}_{n\times p}(\R) \,|\,U^TU=\mathbb{I}_p\}.$$
Denote with  ${\bf u}_1,...,{\bf u}_p\in \R^n$ the  vectors that form the columns of the matrix $U\in  \mathcal{M}_{n\times p}(\R)$. The condition that the matrix $U$ belongs to the orthogonal Stiefel manifold is equivalent with the vectors ${\bf u}_1,...,{\bf u}_p\in \R^n$ being orthonormal.

The functions that describe the constraints defining the orthogonal Stiefel manifold as a preimage of a regular value are $F_{aa},F_{bc}:\mathcal{M}_{n\times p}({\R})\rightarrow \R$ given by:
\begin{align}
F_{aa} (U) & =\frac{1}{2}\|{\bf u}_a\|^2,\,\,1 \leq a\leq p,\nonumber \\
F_{bc} (U) & = \left<{\bf u}_b,{\bf u}_c\right>,\,\,1\leq b<c\leq p. \nonumber 
\end{align}
More precisely, we have ${\bf F}:\mathcal{M}_{n\times p}({\R})\rightarrow \R^{\frac{p(p+1)}{2}}$, ${\bf F}:=\left( \dots , F_{aa},\dots ,F_{bc}, \dots\right)$, $$St_p^n\simeq {\bf F}^{-1}\left( \dots , \frac{1}{2},\dots ,0, \dots\right)\subset \mathcal{M}_{n\times p}({\R}).$$

In what follows we denote by ${\bf vec}(U)\subset \R^{np}$ the column vectorization of the matrix $U\in \mathcal{M}_{n\times p}(\R)$. Also, for a smooth function $f:\mathcal{M}_{n\times p}(\R)\rightarrow \R$ we denote 
$\nabla f(U):=\nabla (f\circ {\bf vec}^{-1})({\bf vec}(U))\in \R^{np}$
and 
$\text{Hess}\,f(U):=\text{Hess}\,(f\circ {\bf vec}^{-1})({\bf vec}(U))\in \mathcal{M}_{np\times np}(\R).$

The Lagrange multiplier functions for the case of an orthogonal Stiefel manifold, see \cite{birtea-casu-comanescu}, are given by the formulas: 
\begin{equation}\label{Lagrange-multipliers-functions}
\left.\begin{array}{l}
\sigma_{aa}(U)=  \left<\nabla G(U),\nabla F_{aa}(U)\right>=\left<\displaystyle\frac{\partial G}{\partial {\bf u}_a}(U),{\bf u}_a\right>;\\ \\
\sigma_{bc}(U)= \left<\nabla G(U),\nabla F_{bc}(U)\right>=\displaystyle\frac{1}{2}\left(\left<\displaystyle\frac{\partial G}{\partial {\bf u}_c}(U),{\bf u}_b\right>+\left<\displaystyle\frac{\partial G}{\partial {\bf u}_b}(U),{\bf u}_c\right>\right),\end{array}\right.
\end{equation}
where $G:\mathcal{M}_{n\times p}({\R})\rightarrow \R$ is the extension of the cost function $\widetilde{G}:St^n_p\rightarrow \R$.

We introduce the symmetric matrix
$$\Sigma(U): =\left[\sigma_{bc}(U)\right]\in \mathcal{M}_{p\times p}(\R),$$
where we define $\sigma_{cb}(U):=\sigma_{bc}(U)$ for $1\leq b<c\leq p$.
Using \eqref{Lagrange-multipliers-functions}, we have\footnote{When we regard $\nabla G(U)$ as a matrix we mean ${\bf vec}^{-1}(\nabla G(U))$.}
\begin{equation}\label{Sigma-matriceal}
\Sigma(U)=\frac{1}{2}\left(\nabla G(U)^TU+U^T\nabla G(U)\right).
\end{equation}

The matrix form of the embedded gradient vector field is given by (see  \cite{birtea-casu-comanescu})
\begin{equation*}
\partial G(U)=\nabla G(U)-U\Sigma(U).
\end{equation*}

The following result provides necessary and sufficient conditions for critical points of a cost function defined on a Stiefel manifold.

\begin{thm}[\cite{birtea-casu-comanescu}]\label{second-egregium}
A matrix $U\in St^n_p$ is a critical point for the cost function $\widetilde{G}=G_{|_{St^n_p}}$ { if and only if} the following conditions are satisfied:
\begin{equation}\label{conditii-matriceale}
\left.
\begin{array}{l}
(i)\,\,U^T\nabla G(U)=\nabla G(U)^TU; \\

(ii)\,\,\nabla G(U)=UU^T\nabla G(U).
\end{array}\right.
\end{equation}
\end{thm}

The Hessian matrices of the constraint functions are given by\footnote{ By $X\otimes Y$ we denote the Kronecker product of the matrices $X,Y$. \\
The vectors ${\bf f}_1$, ...   ,${\bf f}_p$ form the canonical basis in the Euclidean space $\R^p$. The $p\times p$ matrix ${\bf f}_a\otimes {\bf f}^T_b$ has $1$ on the $a$-th row and $b$-th column and the rest $0$. }:
\begin{equation*}
\left[\text{Hess}\,{F_{aa}}(U)\right] =
{\small \begin{blockarray}{cccccc}
&&\left.\begin{array}{c}a~\hbox{col.}\\ \downarrow\end{array}\right.&&\\
\begin{block}{[ccccc]c}
  \mathbb{O}_n & \dots & \mathbb{O}_n & \dots & \mathbb{O}_n &\\
  \dots & \dots & \dots & \dots & \dots &\\
   \mathbb{O}_n & \dots & \mathbb{I}_n & \dots & \mathbb{O}_n &\leftarrow a~\hbox{row}\\
    \dots & \dots & \dots & \dots & \dots &\\
\mathbb{O}_n & \dots & \mathbb{O}_n & \dots & \mathbb{O}_n&\\ \end{block}
\end{blockarray}}
= \left({\bf f}_a\otimes {\bf f}^T_a\right)\otimes \mathbb{I}_n;
\end{equation*}
\begin{align*}
\left[\text{Hess}\,{F_{bc}}(U)\right] & =
{\small \begin{blockarray}{cccccccc}
&&\left.\begin{array}{c}b~\hbox{col.}\\ \downarrow\end{array}\right.&&\left.\begin{array}{c}c~\hbox{col.}\\ \downarrow\end{array}\right.&&\\
\begin{block}{[ccccccc]c}
  \mathbb{O}_n & \dots & \mathbb{O}_n & \dots & \mathbb{O}_n & \dots & \mathbb{O}_n& \\
  \dots & \dots & \dots & \dots & \dots & \dots & \dots&\\
  \mathbb{O}_n &   \dots &  \mathbb{O}_n &  \dots &   \mathbb{I}_n &  \dots &  \mathbb{O}_n&\leftarrow b~\hbox{row}\\
   \dots & \dots & \dots & \dots & \dots & \dots & \dots&\\
  \mathbb{O}_n &   \dots &  \mathbb{I}_n &  \dots &  \mathbb{O}_n  &   \dots &  \mathbb{O}_n&\leftarrow c~\hbox{row}\\
   \dots & \dots & \dots & \dots & \dots & \dots & \dots\\
 \mathbb{O}_n & \dots & \mathbb{O}_n & \dots & \mathbb{O}_n & \dots & \mathbb{O}_n&\\
 \end{block}
\end{blockarray}}\nonumber \\
& = \left({\bf f}_b\otimes {\bf f}^T_c+{\bf f}_c\otimes {\bf f}^T_b\right)\otimes \mathbb{I}_n.
\end{align*}

For the case of orthogonal Stiefel manifold, the general formula for the Hessian matrix of the cost function $\widetilde{G}$, as given in \cite{Birtea-Comanescu-Hessian}, becomes:
$$\text{Hess} \,\widetilde{G}(U):T_U St_p^n\times T_U St_p^n\rightarrow \R,$$

{ \begin{equation}\label{formula-Hessian-234}
 \text{Hess} \widetilde{G}(U)=\left(\text{Hess}G(U)-\sum_{a=1}^p\sigma_{aa}(U)\text{Hess}F_{aa}(U)-\sum_{1\leq b<c\leq p}\sigma_{bc}(U)\text{Hess}F_{bc}(U)\right)_{|T_U St_p^n\times T_U St_p^n}.
\end{equation}}

Using the above expressions for the Hessian matrices of the constraint functions and substituting them in \eqref{formula-Hessian-234} we obtain the following concise matrix form for the Hessian of the cost function.
\begin{thm}\label{Hessian-matricial-123}
The Hessian of the cost function $\widetilde{G}:St_p^n\rightarrow \R$ is given by
\begin{equation*}\label{}
\emph{Hess}\,\widetilde{G}(U) =\left(\emph{Hess}\,G(U)-{\Sigma}(U)\otimes \mathbb{I}_n\right)_{|T_U St_p^n\times T_U St_p^n}.
\end{equation*}
\end{thm}

\section{Local frames on Stiefel manifolds}\label{section-3}
There exist two frequently used methods to prove that a certain set has a manifold structure. One of them is to prove that the desired set is the preimage of a regular value of a smooth function. Another possibility is to explicitly construct compatible local coordinates (local charts)  that cover the entire set. The first approach gives an implicit description of the tangent space. The second approach gives an explicit formula for a basis of the tangent space. 
Regarding the orthogonal Stiefel manifold as a preimage of a regular value, we explicitly construct a local frame on this manifold, although the manifold is described implicitly.
\medskip

\noindent {\bf I. Construction of an explicit local frame for orthogonal Stiefel manifolds.} 

Using the description of a Stiefel manifold as a preimage of a regular value, 
in \cite{edelman} it is given the following elegant explicit form for the tangent space at a point $U\in St_p^n$:
$$T_U St_p^n=\{UA+(\mathbb{I}_n-UU^T)C\,|\,A\in \mathcal{M}_{p\times p}(\R),\,A=-A^T,\,C\in \mathcal{M}_{n\times p}(\R)\}.$$
On the tangent space we consider the Frobenius scalar product:
$$\left<\Delta_1(U),\Delta_2(U)\right>=\tr(\Delta_1^T(U)\Delta_2(U)),\,\,\,\Delta_1(U),\Delta_2(U)\in T_U St_p^n.$$
We construct a basis $\mathcal{B}_U$ for the tangent space $T_U St_p^n$. The explicit description of the vectors in this basis, as we will see below, decisively depends on the choice of a full rank $p\times p$ submatrix of $U$.    We split the basis $\mathcal{B}_U$ as the following set union of vectors
$$\mathcal{B}_U=\mathcal{B}_U^{'}\cup \mathcal{B}_U^{''}.$$
The set $\mathcal{B}_U^{'}$ is formed with tangent vectors of the form 
\begin{equation}\label{delta-prim}
\Delta'_{ab}(U)=UA_{ab},\,\,\,1\leq a<b\leq p,
\end{equation}
where 
$$A_{ab}=(-1)^{a+b}({\bf f}_a\otimes {\bf f}^T_b-{\bf f}_b\otimes {\bf f}_a^T),\,\,\,1\leq a<b\leq p,$$
form the standard basis for the $p\times p$ skew-symmetric matrices. 

For the next computations we use the following rule for matrix multiplication: 
$$\left({\bf u}\otimes {\bf v}^T_{\oplus}\right)\cdot \left({\bf v}_{\odot}\otimes {\bf w}^T\right)=\left({\bf v}^T_{\oplus}\cdot{\bf v}_{\odot}\right) \left({\bf u}\otimes {\bf w}^T\right),$$
where the vectors ${\bf v}_{\oplus}$ and ${\bf v}_{\odot}$ belong to the same space. The matrix product ${\bf v}^T_{\oplus}\cdot{\bf v}_{\odot}$ represents, in fact, the inner product of the vectors ${\bf v}_{\oplus}$ and ${\bf v}_{\odot}$ in the corresponding space. If these two vectors belong to an orthonormal set of vectors, then
$$\left({\bf u}\otimes {\bf v}^T_{\oplus}\right)\cdot \left({\bf v}_{\odot}\otimes {\bf w}^T\right)=\delta_{\oplus,\odot}\left({\bf u}\otimes {\bf w}^T\right),$$
where $\delta_{\oplus,\odot}$ represents the Kronecker delta function. 

\begin{prop}
The tangent vectors in $\mathcal{B}_U^{'}$ are nonzero and orthogonal one to another.
\end{prop}
{\it Proof}
Assume $\Delta'_{ab}(U)=\mathbb{O}_{p\times p}$. Multiplying the equality to the left with the matrix $U^T$, we obtain that $A_{ab}=\mathbb{O}_{p\times p}$ which is a contradiction.

By a direct computation, for $(a_1,b_1)\neq (a_2,b_2)$, we obtain:
\begin{align*}
\left<\Delta'_{a_1b_1}(U),\Delta'_{a_2b_2}(U)\right>  & = \left<UA_{a_1b_1},UA_{a_2b_2}\right> = \tr(A_{a_1b_1}^TU^TUA_{a_2b_2}) = \tr(A_{a_1b_1}^TA_{a_2b_2}) \\
& =\varepsilon\tr(({\bf f}_{b_1}\otimes {\bf f}^T_{a_1}-{\bf f}_{a_1}\otimes {\bf f}^T_{b_1})\cdot ({\bf f}_{a_2}\otimes {\bf f}^T_{b_2}-{\bf f}_{b_2}\otimes {\bf f}^T_{a_2})) \\
& = \varepsilon( \delta_{a_1a_2}\tr({\bf f}_{b_1}\otimes {\bf f}^T_{b_2})-\delta_{a_1b_2}\tr({\bf f}_{b_1}\otimes {\bf f}^T_{a_2}) \\
& -\delta_{a_2b_1}\tr({\bf f}_{a_1}\otimes {\bf f}^T_{b_2})+\delta_{b_1b_2}\tr({\bf f}_{a_1}\otimes {\bf f}^T_{a_2}))= 0,
\end{align*}
where $\varepsilon=(-1)^{a_1+b_1+a_2+b_2}$.\,\rule{0.5em}{0.5em}

As a consequence, we obtain that the set $\mathcal{B}_U^{'}$ has $\frac{p(p-1)}{2}$ linearly independent tangent vectors.

A matrix $U\in St_p^n$ has rank $p$ and the construction of tangent vectors in $\mathcal{B}_U^{''}$ crucially depends on the choice of a $p\times p$ full rank submatrix of $U$. Assume that this $p\times p$ full rank submatrix $U_p$ is formed with the rows $1\leq i_1<i_2<...<i_p\leq n$ of the matrix $U$.

Define the tangent vectors of $\mathcal{B}_U^{''}$ as
\begin{equation}\label{delta-secund}
\Delta''_{ic}(U)=(\mathbb{I}_n-UU^T)C_{ic},\,\,i\in \{1,...,n\}\backslash\{i_1,...,i_p\},\,c\in \{1,...,p\},
\end{equation}
where 
$$C_{ic}={\bf e}_i\otimes {\bf f}^T_c.\footnote{The vectors ${\bf e}_1$, ...   ,${\bf e}_n$ form the canonical basis in the Euclidean space $\R^n$. The $n\times p$ matrix ${\bf e}_i\otimes {\bf f}^T_c$ has $1$ on the $i$-th row and $c$-th column and the rest $0$. }$$
For a fixed $c\in \{1,...,p\}$, we define:
\begin{equation}\label{secund-3}
 _c\mathcal{B}_{U}^{''}=\{\Delta''_{ic}(U)\,|\,i\in \{1,...,n\}\backslash\{i_1,...,i_p\}\}.
\end{equation}
The set $\mathcal{B}_U^{''}$ is defined by 
\begin{equation}\label{secund-2}
\mathcal{B}_U^{''}=\bigcup_{c=1}^{p}{_c\mathcal{B}_{U}^{''}}.
\end{equation}

\begin{prop} The tangent vectors in $\mathcal{B}_{U}^{''}$ have the following properties:
\begin{itemize}
\item [(i)] For a fixed $c\in \{1,...,p\}$ the tangent vectors in $_c\mathcal{B}_{U}^{''}$ are linearly independent.
\item [(ii)] For $c_1,c_2\in \{1,...,p\}$ and $c_1\neq c_2$ we have 
$$\left<\Delta''_{k_1c_1}(U),\Delta''_{k_2c_2}(U)\right>=0,\,\,\,\forall\, k_1,k_2\in \{1,...,n\}\backslash\{i_1,...,i_p\}.$$
\end{itemize} 
\end{prop}
{\it Proof} $(i)$ We make the notation $I_p=\{i_1,...,i_p\}$. Using the tensorial description of a $n\times p$ matrix we have:
\begin{align*}
U= & \sum_
{\substack{
   j\notin I_p \\
b\in \{1,...,p\}
  }}
u_{jb}\,{\bf e}_j\otimes {\bf f}^T_b+
\sum_
{\substack{
   k\in I_p \\
a\in \{1,...,p\}
  }}
u_{ka}\,{\bf e}_k\otimes {\bf f}^T_a, 
\end{align*}

\begin{align*}
U^T= & \sum_
{\substack{
   s\notin I_p \\
d\in \{1,...,p\}
  }}
u_{sd}\,{\bf f}_d\otimes {\bf e}^T_s+
\sum_
{\substack{
   r\in I_p \\
g\in \{1,...,p\}
  }}
u_{rg}\,{\bf f}_g\otimes {\bf e}^T_r, 
\end{align*}
where ${\bf e}_j,{\bf e}_k,{\bf e}_r,{\bf e}_s\in \R^n$ and ${\bf f}_a,{\bf f}_b,{\bf f}_d,{\bf f}_g\in \R^p$. 
 By a direct computation we  obtain the tensorial description of $UU^T$:
\begin{align*}
UU^T= & \sum_
{\substack{
   j,s\notin I_p \\
b\in \{1,...,p\}
  }}
u_{jb}u_{sb}\,{\bf e}_j\otimes {\bf e}^T_s+
\sum_
{\substack{
   j\notin I_p,\,r\in I_p \\
b\in \{1,...,p\}
  }}
u_{jb}u_{rb}\,{\bf e}_j\otimes {\bf e}^T_r \\
& +  \sum_
{\substack{
   k\in I_p,\,s\notin I_p \\
a\in \{1,...,p\}
  }}
u_{ka}u_{sa}\,{\bf e}_k\otimes {\bf e}^T_s+
\sum_
{\substack{
   k,r\in I_p \\
a\in \{1,...,p\}
  }}
u_{ka}u_{ra}\,{\bf e}_k\otimes {\bf e}^T_r.
\end{align*}
For proving the linear independence of the vectors in $_c\mathcal{B}_{U}^{''}$, we have the computation:
\begin{align*}
\mathbb{O}_{n\times p}  = &  \sum_{i\notin I_p}\alpha_i\Delta''_{ic}(U)=\sum_{i\notin I_p}\alpha_i(\mathbb{I}_n-UU^T)C_{ic} \\
  = &  \sum_{i\notin I_p}\alpha_i{\bf e}_i\otimes {\bf f}^T_c
-\sum_
{\substack{
   i,j,s\notin I_p \\
b\in \{1,...,p\}
  }}
\alpha_i u_{jb}u_{sb}\,\delta_{is}{\bf e}_j\otimes {\bf f}^T_c
-\sum_
{\substack{
   i,j\notin I_p,\,r\in I_p \\
b\in \{1,...,p\}
  }}
\alpha_i u_{jb}u_{rb}\,\delta_{ir}{\bf e}_j\otimes {\bf f}^T_c \\
& - \sum_
{\substack{
   i,s\notin I_p,\,k\in I_p \\
a\in \{1,...,p\}
  }}
\alpha_i u_{ka}u_{sa}\,\delta_{is}{\bf e}_k\otimes {\bf f}^T_c
 - \sum_
{\substack{
   i\notin I_p,\,k,r\in I_p \\
a\in \{1,...,p\}
  }}
\alpha_i u_{ka}u_{ra}\,\delta_{ir}{\bf e}_k\otimes {\bf f}^T_c \\
=& \sum_{i\notin I_p}\alpha_i{\bf e}_i\otimes {\bf f}^T_c 
-\sum_
{\substack{
   i,j\notin I_p \\
b\in \{1,...,p\}
  }}
\alpha_i u_{jb}u_{ib}\,{\bf e}_j\otimes {\bf f}^T_c
-\sum_
{\substack{
   i\notin I_p,\,k\in I_p \\
a\in \{1,...,p\}
  }}
\alpha_i u_{ka}u_{ia}\,{\bf e}_k\otimes {\bf f}^T_c \\
= &  \sum_{j\notin I_p}\left( \alpha_j-\sum_
{\substack{
   i\notin I_p \\
b\in \{1,...,p\}
  }}
\alpha_i u_{jb}u_{ib}\right) \,{\bf e}_j\otimes {\bf f}^T_c -
\sum_{k\in I_p}\left( \sum_
{\substack{
   i\notin I_p \\
a\in \{1,...,p\}
  }}
\alpha_i u_{ka}u_{ia}\right) {\bf e}_k\otimes {\bf f}^T_c,
\end{align*}
where $\delta_{ir}=0$ for any $i\notin I_p$ and $r\in I_p$. Decomposing the above matrix equality on the subspaces $\text{Span} \{{\bf e}_j\otimes {\bf f}^T_c\,|\,j\notin I_p\}$ and  $\text{Span} \{{\bf e}_k\otimes {\bf f}^T_c\,|\,k\in I_p\}$ we obtain:
\begin{align}
&\alpha_j-\sum\limits_
{\substack{
  i\notin I_p \\
b\in \{1,...,p\}
  }}
\alpha_i u_{jb}u_{ib}=0,\,\,\,\forall j\notin I_p\label{ec-unu} \\
& \sum\limits_
{\substack{
   i\notin I_p \\
a\in \{1,...,p\}
  }}
\alpha_i u_{ka}u_{ia}=0,\,\,\,\forall k\in I_p. \label{ec-doi}
\end{align}
In order to write the above system in a matrix form we need to relabel the elements of the set $\{1,...,n\}\backslash I_p$. There exists a unique strictly increasing function $\sigma: \{1,...,n\}\backslash I_p\rightarrow \{1,...,n-p\}$. Analogously, we relabel the set $I_p$ using the unique strictly increasing function $\tau: I_p\rightarrow \{1,...,p\}$. 

The $p\times p$ full rank submatrix $U_p$ has the following tensorial form:
$$U_p=\sum\limits_
{\substack{
   k\in I_p \\
a\in \{1,...,p\}
  }}
u_{ka}\,{\bf f}_{\tau(k)}\otimes {\bf f}^T_a$$
and the $(n-p)\times p$ complement $U_{n-p}$ of the full rank submatrix $U_p$ in the matrix $U$ has the following tensorial form:
$$U_{n-p}=\sum\limits_
{\substack{
   j\notin I_p \\
a\in \{1,...,p\}
  }}
u_{ja}\,{\bf h}_{\sigma(j)}\otimes {\bf f}^T_a.\footnote{The vectors ${\bf h}_{\sigma(i)}$ form the canonical basis of $\R^{n-p}$. }$$
We define the following $(n-p)\times p$ submatrix of $C_{ic}$ removing the rows $i_1,\dots,i_p$,
$$\widetilde{C}_{ic}={\bf h}_{\sigma(i)}\otimes {\bf f}^T_c.$$
The equations \eqref{ec-unu} and \eqref{ec-doi} have the equivalent forms:
$$\sum\limits_
{\substack{
  j\notin I_p
  }}
\left(\alpha_j\widetilde{C}_{jc}-\alpha_j U_{n-p}U_{n-p}^T\widetilde{C}_{jc}\right)=\mathbb{O}_{(n-p)\times p}
$$
and respectively
$$\sum\limits_
{\substack{
  i\notin I_p
  }}
\alpha_iU_pU_{n-p}^T\widetilde{C}_{ic}=\mathbb{O}_{p\times p}.
$$
Because the full rank submatrix $U_p$ is invertible we obtain
$\sum\limits_
{\substack{
  i\notin I_p
  }}
\alpha_iU_{n-p}^T\widetilde{C}_{ic}=\mathbb{O}_{p\times p}$ and consequently, 
$\sum\limits_
{\substack{
  i\notin I_p
  }}
\alpha_i\widetilde{C}_{ic}=\mathbb{O}_{(n-p)\times p}$. The linear independence of the matrices $\widetilde{C}_{ic}$ implies $\alpha_i=0$ for all $i\notin I_p$, which proves the linear independence of the vectors in $_c\mathcal{B}_{U}^{''}$.

$(ii)$ We make the notation $Z=\mathbb{I}_n-UU^T$ and a direct computation shows that $Z^TZ=Z$. Consequently, when $c_1\neq c_2$ we have:
\begin{align*}
\left<\Delta''_{k_1c_1}(U),\Delta''_{k_2c_2}(U)\right>= & \tr (C_{k_1c_1}^TZ^TZC_{k_2c_2})= \tr (C_{k_1c_1}^TZC_{k_2c_2}) \\
= & \tr\left(({\bf f}_{c_1}\otimes {\bf e}^T_{k_1})\cdot (\sum_{i,j}z_{ij}{\bf e}_{i}\otimes {\bf e}^T_{j})\cdot ({\bf e}_{k_2}\otimes {\bf f}^T_{c_2})\right) \\
= & \tr\left(z_{k_1k_2}{\bf f}_{c_1}\otimes {\bf f}^T_{c_2}\right)=0.\,\,\,\rule{0.5em}{0.5em}
\end{align*}

\begin{prop}
The tangent vectors in $\mathcal{B}_{U}^{'}$ are orthogonal to the vectors in $\mathcal{B}_{U}^{''}$.
\end{prop}
{\it Proof} $\,$ We notice that $U^T(\mathbb{I}_n-UU^T)=\mathbb{O}_{p\times n}$. Consequently, for $1\leq a<b\leq p$, $i\notin I_p$ and $c\in \{1,...,p\}$ we have:
\begin{align*}
\left<\Delta'_{ab}(U),\Delta''_{ic}(U)\right>= & \tr (A_{ab}^TU^T(\mathbb{I}_n-UU^T)C_{ic})=0.\,\,\,\rule{0.5em}{0.5em}
\end{align*}

The above results can be summarized in the following theorem.

\begin{thm}
For $U\in St^n_p$ the vectors of the set $\mathcal{B}_U=\mathcal{B}_U^{'}\cup \mathcal{B}_U^{''}$, where the elements of the set $\mathcal{B}_U^{'}$ are given by \eqref{delta-prim} and the elements of the set $\mathcal{B}_U^{''}$ are given by 
\eqref{delta-secund}, \eqref{secund-3}, and \eqref{secund-2},  form a basis for the tangent space $T_U St_p^n$. Among them, we have the following orthogonality properties:
\begin{itemize}
\item [(i)] $\mathcal{B}_U^{'}$ is an orthogonal set;
\item [(ii)] for $c_1,c_2\in \{1,...,p\}$ and $c_1\neq c_2$ we have $_{c_1}\mathcal{B}_{U}^{''}\perp {_{c_2}\mathcal{B}_{U}^{''}}$;
\item [(iii)] $\mathcal{B}_U^{'}\perp \mathcal{B}_U^{''}$.
\end{itemize}
\end{thm}

In order to have all vectors in the basis $\mathcal{B}_U$ orthogonal one to the other, we can apply the Gram-Schmidt algorithm to the vectors in $_{c}\mathcal{B}_{U}^{''}$ and do this for each $c\in \{1,...,p\}$.

The above construction of the vectors in $\mathcal{B}_{U}^{''}$ crucially depends on the choice of the full rank $p\times p$ submatrix $U_p$. Choosing another full rank submatrix leads to a change of basis for the tangent space $T_U St_p^n$.
\medskip

\noindent {\bf II. Local frames on sphere $S^{n-1}$}. For the case when $p=1$, the Stiefel manifold $St^n_1$ becomes the sphere $S^{n-1}\subset \R^{n}$.
In this case, for ${\bf x}\in S^{n-1}$, we have $\mathcal{B}'_{\bf x}=\emptyset$. For a point ${\bf x}\in S^{n-1}$, we choose an index $j\in \{1,...,n\}$ such that $x_j\neq 0$. Consequently, a local frame for the sphere is given by 
$$\Delta_{i1}''({\bf x})=(\mathbb{I}_n-{\bf x\otimes x}^T){\bf e}_i\otimes {\bf f}^T_1=(\mathbb{I}_n-{\bf x\otimes x}^T){\bf e}_i={\bf e}_i-x_i{\bf x},\,\,\,i\in \{1,...,n\}\setminus \{j\}.$$
{\bf III. The Hessian of the constraint functions computed on the basis $\mathcal{B}_U$.} The orthogonality among the elements of the basis $\mathcal{B}_U$ has the computational advantage that it renders the Hessian matrices of the constraint functions in a very simple form, where most of the entries are zero.

More precisely, let $\Delta'_{\alpha_1\beta_1}(U)$ and $\Delta'_{\alpha_2\beta_2}(U)$ be two elements of  $\mathcal{B}'_U$. By a direct computation, we have the following formulas:
\begin{align*}
\text{Hess}\, F_{aa}\left(\Delta'_{\alpha_1\beta_1}(U),\Delta'_{\alpha_2\beta_2}(U)\right) = &
(-1)^{\alpha_1+\beta_1+\alpha_2+\beta_2}\left(\delta_{a\alpha_1}\delta_{a\alpha_2}\delta_{\beta_1\beta_2}+\delta_{a\beta_1}\delta_{a\beta_2}\delta_{\alpha_1\alpha_2}\right); \\
\text{Hess}\, F_{bc}\left(\Delta'_{\alpha_1\beta_1}(U),\Delta'_{\alpha_2\beta_2}(U)\right)  = &
(-1)^{\alpha_1+\beta_1+\alpha_2+\beta_2}\left(\delta_{\alpha_1 b}\delta_{c\alpha_2}\delta_{\beta_1\beta_2}-\delta_{\alpha_1 b}\delta_{c\beta_2}\delta_{\beta_1\alpha_2}+\delta_{\alpha_1 c}\delta_{b\alpha_2}\delta_{\beta_1\beta_2}    \right.\nonumber \\ 
 & \left.+ \delta_{\beta_1 b}\delta_{c\beta_2}\delta_{\alpha_1\alpha_2}-\delta_{\beta_1 c}\delta_{b\alpha_2}\delta_{\alpha_1\beta_2}+\delta_{\beta_1 c}\delta_{b\beta_2}\delta_{\alpha_1\alpha_2}\right).
\end{align*} 
Further analyzing the above formulas we have:
\begin{equation*}
\text{Hess}\, F_{aa}\left(\Delta'_{\alpha_1\beta_1}(U),\Delta'_{\alpha_2\beta_2}(U)\right) =
\begin{cases}
0\,\,\,\text{if}\,(\alpha_1,\beta_1)\neq (\alpha_2,\beta_2) \\
0\,\,\,\text{if}\,(\alpha_1,\beta_1)= (\alpha_2,\beta_2),\,a\notin \{\alpha_1, \beta_1\} \\
1\,\,\,\text{if}\,(\alpha_1,\beta_1)= (\alpha_2,\beta_2),\,a\in \{\alpha_1, \beta_1\}.
\end{cases}
\end{equation*}
Also, the value $\text{Hess}\, F_{bc}\left(\Delta'_{\alpha_1\beta_1}(U),\Delta'_{\alpha_2\beta_2}(U)\right)$ is 0 or $\pm1$ depending on the ordering and relative position of the integer numbers $b,c,\alpha_1,\beta_1, \alpha_2,\beta_2$.

For the case when $\Delta'_{\alpha\beta}(U)\in \mathcal{B}'(U)$ and $\Delta''_{j d}(U)\in \mathcal{B}''(U)$ we have: 
\begin{equation*}
\text{Hess}\, F_{aa}\left(\Delta'_{\alpha\beta}(U),\Delta''_{j d}(U)\right)=\text{Hess}\, F_{bc}\left(\Delta'_{\alpha\beta}(U),\Delta''_{j d}(U)\right)=0.
\end{equation*}

Let $\Delta''_{j_1 d_1}(U)$ and $\Delta''_{j_2 d_2}(U)$ be two elements of  $\mathcal{B}''_U$. By a direct computation, using the notation $Z=[z_{kl}]=\mathbb{I}_n-UU^T$, we have the following formulas:
\begin{align*}
\text{Hess}\, F_{aa}\left(\Delta''_{j_1 d_1}(U),\Delta''_{j_2 d_2}(U)\right) = &
\delta_{a d_1}\delta_{a d_2}z_{j_1 j_2}=\begin{cases}
z_{j_1j_2}\,\,\,\text{if}\,\,\,a=d_1=d_2 \\
0\,\,\,\,\,\,\,\,\,\,\,\,\text{otherwise}
\end{cases};\\
\text{Hess}\, F_{bc}\left(\Delta''_{j_1 d_1}(U),\Delta''_{j_2 d_2}(U)\right) = &
\left(\delta_{b d_1}\delta_{c d_2}+\delta_{b d_2}\delta_{c d_1}\right) z_{j_1 j_2}=
\begin{cases}
z_{j_1j_2}\,\text{if}\,\,(b,c)=(d_1,d_2)\,\text{or}\, (b,c)=(d_2,d_1)\\
0\,\,\,\,\,\,\,\,\,\,\,\text{otherwise}.
\end{cases}
\end{align*}

Taking into account the matrix form for the Hessian of the cost function from Theorem \ref{Hessian-matricial-123} we could also be interested in the values of the term $\Sigma(U)\otimes \mathbb{I}_n$ on two tangent  vectors from $T_USt^n_p$. These computations are presented in what follows.

For two tangent vectors $\Delta_1(U)=UA_1+(\mathbb{I}_n-UU^T)C_1,\Delta_2(U)=UA_2+(\mathbb{I}_n-UU^T)C_2\in T_USt^n_p$ we have
$$(\Sigma(U)\otimes \mathbb{I}_n)(\Delta_1(U),\Delta_2(U))=-\tr(A_1A_2\Sigma)+\tr(C_1^TZC_2\Sigma).$$

Let $\Delta'_{\alpha\beta}(U)$ be an element of $\mathcal{B}'_U$ and $\Delta''_{j d}(U)$ be an element of $\mathcal{B}''_U$. Then 
\begin{equation}\label{sigma-I-0}
(\Sigma(U)\otimes \mathbb{I}_n)(\Delta'_{\alpha\beta}(U),\Delta''_{j d}(U))=0.
\end{equation}
Let $\Delta'_{\alpha_1\beta_1}(U)$, $\Delta'_{\alpha_2\beta_2}(U)$ be two elements of  $\mathcal{B}'_U$ and $\Delta''_{j_1 d_1}(U)$, $\Delta''_{j_2 d_2}(U)$ be two elements of  $\mathcal{B}''_U$. Then we have
\begin{align}
(\Sigma(U)\otimes \mathbb{I}_n)(\Delta'_{\alpha_1\beta_1}(U),\Delta'_{\alpha_2\beta_2}(U)) & =(-1)^{\alpha_1+\beta_1+\alpha_2+\beta_2}\cdot \left(\delta_{\alpha_1\alpha_2}\sigma_{\beta_1\beta_2}+\delta_{\beta_1\beta_2}\sigma_{\alpha_1\alpha_2}\right.\nonumber\\
& \left.-\delta_{\alpha_2\beta_1}\sigma_{\alpha_1\beta_2}-\delta_{\alpha_1\beta_2}\sigma_{\alpha_2\beta_1}\right)\label{sigma-I-1}
\end{align}
and
\begin{equation}
(\Sigma(U)\otimes \mathbb{I}_n)(\Delta''_{j_1 d_1}(U),\Delta''_{j_2 d_2}(U))=z_{j_1j_2}\sigma_{d_1d_2}.\label{sigma-I-2}
\end{equation}

\section{Applications}

\subsection{The Procrustes problem on orthogonal Stiefel manifolds}

The optimization problem is the following:
$$\underset{\mathlarger{U\in St_p^n}}{\text{Minimize}}\,\,\, ||AU-B||^2,$$
where $A\in \mathcal{M}_{m\times n}(\R)$, $B\in \mathcal{M}_{m\times p}(\R)$,  and $\|\cdot \|$ is the Frobenius norm. The cost function associated to this optimization problem is given by $\widetilde{G}:St^n_p\rightarrow \R$ and its natural extension $G:\mathcal{M}_{n\times p}(\R)\rightarrow \R$  is
$$G(U)=\frac{1}{2}||AU-B||^2=\frac{1}{2}\hbox{tr}(U^TA^TAU)-\hbox{tr}(U^TA^TB)+\frac{1}{2}\hbox{tr}(B^TB).$$
By a straightforward computation we have that
$\nabla G(U)=A^TAU-A^TB.$

First order optimality necessary and sufficient conditions are given in \cite{chu} and \cite{birtea-casu-comanescu} and can be obtained using Theorem \ref{second-egregium}.
\begin{thm}
A matrix $U\in St^n_p$ is a critical point for the Procrustes cost function if and only if:
\begin{itemize}
\item [(i)] the matrix $B^TAU$ is symmetric;
\item [(ii)] $(\mathbb{I}_n-UU^T)(A^TAU-A^TB)={\bf 0}$.
\end{itemize}
\end{thm}

\noindent The following formula holds:
$\text{Hess}\,G(U)=\mathbb{I}_p\otimes (A^TA)$.

\noindent For a tangent vector $\Delta(U)=U K+(\mathbb{I}_n-UU^T) W\in T_USt^n_p$, 
we have the following computation\footnote{${\bf vec}(A)^T(D\otimes B){\bf vec}(C)=\tr(A^TBCD^T)$.}:
\begin{align*}
\text{Hess} \,\widetilde{G}(U)(\Delta(U), \Delta(U)) = & 
  (\mathbb{I}_p\otimes (A^TA)-\Sigma(U)\otimes \mathbb{I}_n)(U K+(\mathbb{I}_n-UU^T)W, UK+(\mathbb{I}_n-UU^T)W) \\
  = & -\tr(KU^TA^TAUK)-2\tr(KU^TA^TA(\mathbb{I}_n-UU^T)W) \\
  & + \tr(W^T(\mathbb{I}_n-UU^T)A^TA(\mathbb{I}_n-UU^T)W) \\
 & +\tr(K^2\Sigma(U))-\tr(W^T(\mathbb{I}_n-UU^T)W\Sigma(U)).
\end{align*}
Using the equality (3.8) from \cite{birtea-casu-comanescu} (see also equation \eqref{Sigma-matriceal}) and the condition $(i)$ of Theorem 3.3 from \cite{birtea-casu-comanescu}, the Lagrange multipliers matrix in a critical point $U$ is given by  
\begin{equation*}
\Sigma(U)=U^TA^TAU-U^TA^TB.
\end{equation*}
Consequently, a necessary condition for a critical point $U\in St^n_p$ to be local minimum is
\begin{align*}
& \left<B^TAUK,K\right>+2\left<A^TAUK,(\mathbb{I}_n-UU^T)W\right>+\left<A^TA(\mathbb{I}_n-UU^T)W,(\mathbb{I}_n-UU^T)W\right> \\
& -\left<(\mathbb{I}_n-UU^T)WU^TA^T(AU-B),(\mathbb{I}_n-UU^T)W\right>\geq 0
\end{align*} 
for all skew-symmetric matrices $K\in \mathcal{M}_{p\times p}(\R)$ and all matrices $W\in \mathcal{M}_{n\times p}(\R)$.

A sufficient condition for a critical point $U\in St^n_p$ to be local minimum is that the above condition is a  strict inequality  for all skew-symmetric matrices $K\in \mathcal{M}_{p\times p}(\R)$ and all matrices $W\in \mathcal{M}_{n\times p}(\R)$,
 such that $K$ and $W$ are not simultaneously null matrices. 
 
 This condition has been previously presented in Theorem 6 from \cite{chu} using a different method. Necessary and sufficient conditions for local and global minimum for the Procrustes problem are obtained in \cite{elden-park}, using the classical Lagrange multipliers method.

\subsection{The Penrose regression problem on orthogonal Stiefel manifolds}

The optimization problem is the following:
$$\underset{\mathlarger{U\in St_p^n}}{\text{Minimize}}\,\,\,  ||AUC-B||^2,$$
where $A\in \mathcal{M}_{m\times n}(\R)$, $B\in \mathcal{M}_{m\times q}(\R)$, $C\in \mathcal{M}_{p\times q}(\R)$.  The cost function is given by $\widetilde{G}:St^n_p\rightarrow \R$ and its natural extension $G:\mathcal{M}_{n\times p}(\R)\rightarrow \R$ is
$$G(U)=\frac{1}{2}||AUC-B||^2=\frac{1}{2}\hbox{tr}(C^TU^TA^TAUC)-\hbox{tr}(C^TU^TA^TB)+\frac{1}{2}\hbox{tr}(B^TB).$$
We have that (see \cite{birtea-casu-comanescu})
$$\nabla G(U)=A^T(AUC-B)C^T.$$ 
The necessary and sufficient conditions of Theorem \ref{second-egregium} for critical points become:

\begin{thm}[\cite{chu}, \cite{birtea-casu-comanescu}]
A matrix $U\in St^n_p$ is a critical point for the Penrose regression cost function if and only if:
\begin{itemize}
\item [(i)] the matrix $C(AUC-B)^TAU$ is symmetric;
\item [(ii)] $(\mathbb{I}_n-UU^T)A^T(AUC-B)C^T={\bf 0} .$
\end{itemize}
\end{thm}

The following formula holds:
$\text{Hess}\,G(U)=(CC^T)\otimes (A^TA)$.

Using the equality (3.8) from \cite{birtea-casu-comanescu} (see also equation \eqref{Sigma-matriceal}) and the condition $(i)$ of Theorem 3.3 from \cite{birtea-casu-comanescu}, the Lagrange multipliers matrix in a critical point $U$ is given by  
\begin{equation*}
\Sigma(U)=U^TA^T(AUC-B)C^T.
\end{equation*}

Using Theorem \ref{Hessian-matricial-123}, the Hessian of the Penrose regression cost function $\widetilde{G}:St_p^n\rightarrow \R$ is given by
\begin{equation*}\label{}
\text{Hess}\,\widetilde{G}(U) =\left((CC^T)\otimes (A^TA)-(U^TA^T(AUC-B)C^T)\otimes \mathbb{I}_n\right)_{|T_U St_p^n\times T_U St_p^n}.
\end{equation*}
Applying the above formula we can deduce a necessary and sufficient condition  for a critical point $U\in St^n_p$ to be a local minimum. This condition has been previously obtained in Theorem 4 from \cite{chu} using another method.

\subsection{The Brockett problem on orthogonal Stiefel manifolds}

The optimization problem is the following, see \cite{absil-mahony-sepulchre-1}, \cite{birtea-casu-comanescu}:
$$\underset{\mathlarger{U\in St_p^n}}{\text{Minimize}}\,\,\, \tr (U^TAUN),$$
where $A\in \mathcal{M}_{n\times n}(\R)$ is a symmetric matrix, $N=\text{diag}(\mu_1,\dots,\mu_p)\in \mathcal{M}_{p\times p}(\R)$ with $0\leq \mu_1\leq\dots\leq \mu_p$. The cost function is given by $\widetilde{G}:St^n_p\rightarrow \R$ and its natural extension $G:\mathcal{M}_{n\times p}(\R)\rightarrow \R$ is 
$$G(U)= \tr (U^TAUN).$$
 By direct computation, see also \cite{absil-mahony-sepulchre-1},  we obtain $\nabla G(U)=2AUN$.
 
 The necessary and sufficient conditions of Theorem \ref{second-egregium} for critical points become:

\begin{thm}[\cite{absil-mahony-sepulchre-1}, \cite{birtea-casu-comanescu}]
A matrix $U\in St^n_p$ is a critical point for the Brockett  cost function if and only if:
\begin{itemize}
\item [(i)] the matrix $U^TAUN$ is symmetric;
\item [(ii)] $AUN=UU^TAUN .$
\end{itemize}
\end{thm}

The following formula holds:
$\text{Hess}\,G(U)=2N\otimes A$.

The Lagrange multipliers matrix in a critical point $U$ is given by  
\begin{equation*}
\Sigma(U)=2U^TAUN.
\end{equation*}

Using Theorem \ref{Hessian-matricial-123}, the Hessian of the Brockett cost function $\widetilde{G}:St_p^n\rightarrow \R$ is given by
\begin{equation}\label{Hessian-Brockett}
\text{Hess}\,\widetilde{G}(U) =2\left(N\otimes A-(U^TAUN)\otimes \mathbb{I}_n\right)_{|T_U St_p^n\times T_U St_p^n}.
\end{equation}

By similar computations as for the Procrustes problem we obtain necessary and sufficient conditions for a critical point $U \in St_p^n$ to be a local minimum for the Brockett cost function. More precisely,  
if the critical point $U \in St_p^n$ is a local minimum, then 
\begin{align*}
& \tr(W^T(\mathbb{I}_n-UU^T)A(\mathbb{I}_n-UU^T)WN)-\tr(KU^TAUKN) \\
& -2\tr(KU^TA(\mathbb{I}_n-UU^T)WN)
-\tr((-K^2+W^T(\mathbb{I}_n-UU^T)W)U^TAUN)\geq 0
\end{align*}
for all skew-symmetric matrices $K\in \mathcal{M}_{p\times p}(\R)$ and all matrices $W\in \mathcal{M}_{n\times p}(\R)$.

A sufficient condition for a critical point $U\in St^n_p$ to be local minimum is that the above condition is a  strict inequality  for all skew-symmetric matrices $K\in \mathcal{M}_{p\times p}(\R)$ and all matrices $W\in \mathcal{M}_{n\times p}(\R)$,
 such that $K$ and $W$ are not simultaneously null matrices.

In what follows we study the particular  Brockett cost function defined  on $St^4_2$ by $A=\text{diag}\,(1,2,3,4)$ and $N=\text{diag}\,(1,2)$. We give the list of all critical points and we completely characterize them. 

The cost function being quadratic it is invariant under the sign change of the vectors that give the columns of the matrix $U$, but it is not invariant under the order of these column vectors.
As shown in \cite{birtea-casu-comanescu}, a matrix $U\in St^4_2$ is a critical point of the Brockett cost function if and only if every column vector of the matrix $U$ is an eigenvector of the matrix $A$.
Computing the eigenvalues of the Hessian matrix using the formula \eqref{Hessian-Brockett}, we obtain the following characterization of the critical points for the above Brockett cost function:
\begin{itemize}[leftmargin=0.5cm]
\item four global minima generated by $[{\bf e}_2,{\bf e_1}]$ (i.e., $[{\bf e}_2,{\bf e_1}]$, $[-{\bf e}_2,{\bf e_1}]$, $[{\bf e}_2,-{\bf e_1}]$, and $[-{\bf e}_2,-{\bf e_1}]$) with the value of the cost function equals 4.
\item eight saddle points generated by $[{\bf e}_1,{\bf e_2}]$ and $[{\bf e}_3,{\bf e_1}]$ with the value of the cost function equals 5.
\item four saddle points generated by $[{\bf e}_4,{\bf e_1}]$ with the value of the cost function equals 6.
\item eight saddle points generated by $[{\bf e}_1,{\bf e_3}]$ and $[{\bf e}_3,{\bf e_2}]$ with the value of the cost function equals 7.
\item eight saddle points generated by $[{\bf e}_2,{\bf e_3}]$ and $[{\bf e}_4,{\bf e_2}]$ with the value of the cost function equals 8.
\item four saddle points generated by $[{\bf e}_1,{\bf e_4}]$ with the value of the cost function equals 9.
\item eight saddle points generated by $[{\bf e}_2,{\bf e_4}]$ and $[{\bf e}_4,{\bf e_3}]$ with the value of the cost function equals 10.
\item four global maxima generated by $[{\bf e}_3,{\bf e_4}]$ with the value of the cost function equals 11.
\end{itemize}

In \cite{dodig} it is discussed a quadratic programing problem, which is equivalent with the following optimization problem on the Stiefel manifold $St^3_2$: 
$$\underset{\mathlarger{U\in St_2^3}}{\text{Minimize}}\,\,\, {\bf vec}(U)^TC\,{\bf vec}(U),$$
where $C\in \mathcal{M}_{6\times 6}(\R)$. In the particular case when the matrix $C$ can be written in the form $C=C_1\otimes C_2$ with $C_1\in \mathcal{M}_{2\times 2}(\R)$ and $C_2\in \mathcal{M}_{3\times 3}(\R)$, the cost function associated to the optimization problem has the expression $\widetilde{G}(U)=\tr (U^TC_2UC_1^T)$. Notice that this cost function is similar with the Brockett cost function, if we drop the symmetry condition for the matrices $A$ and $N$. The computations for the first and second order optimality discussion can be carried out in an analogous  manner as above.

\subsection{The Newton algorithm on Stiefel manifolds}\label{section-4}

We first recall the setting of the Newton method presented in \cite{5-electron}.
The iterative scheme of Newton algorithm on Riemannian manifolds is given by: 
\begin{equation*}
x_{k+1}=\wt{\mathcal{R}}_{x_k}(\wt{{\bf v}}_{x_k}),
\end{equation*}
where the sequence $\left(x_k\right)_{k\in \mathbb{N}}$ belongs to a smooth Riemannian submanifold $(\mathcal{S},{\bf g}_{_{ind}})$ embedded in a larger manifold $(\mathcal{M},{\bf g})$ (as in the general context presented in Section 2), $\wt{G}:\mathcal{S}\rightarrow \R$ is a smooth cost function, $\wt{\mathcal{R}}:T\mathcal{S}\rightarrow \mathcal{S}$ is a smooth retraction, and the tangent vector $\wt{{\bf v}}_{x_k}\in T_{x_k}\mathcal{S}$ is the solution of the {\bf (contravariant) Newton equation}
\begin{equation}\label{Newton-equation-standard}
\mathcal{H}^{\wt{G}}(x_k)\cdot \wt{{\bf v}}_{x_k}=-\nabla_{{\bf g}_{_{ind}}}\wt{G}(x_k).
\end{equation}
Using the link between the Hessian operator $\mathcal{H}^{\wt{G}}:T\mathcal{S}\rightarrow T\mathcal{S}$ and its associated symmetric bilinear form $\text{Hess}\,{\wt{G}}:T\mathcal{S}\times T\mathcal{S}\rightarrow \R$, the equation \eqref{Newton-equation-standard} can be written equivalently as the {\bf (covariant) Newton equation}: 
\begin{equation}\label{Newton-equation-forms}
 \text{Hess}\, {\wt{G}}(x_k)\cdot \wt{{\bf v}}_{x_k}=-d\wt{G}(x_k).
\end{equation}

In what follows, we customize the Embedded Newton Algorithm from \cite{5-electron} to the specific case of orthogonal Stiefel manifolds.
\medskip

\noindent {\bf Embedded Newton algorithm on Stiefel manifolds}:
\begin{framed}
\begin{itemize}
\item [1.] Consider a smooth prolongation $G:\mathcal{M}_{n\times p}(\R)\rightarrow \R$ of the cost function  $\wt{G}:St^n_p\rightarrow \R$.

\item [2.] Compute the Lagrange multiplier functions:
\begin{equation*}
		\left.\begin{array}{l}
		\sigma_{aa}(U)=  \left<\displaystyle\frac{\partial G}{\partial {\bf u}_a}(U),{\bf u}_a\right>,\,\,
		\sigma_{bc}(U)=\displaystyle\frac{1}{2}\left(\left<\displaystyle\frac{\partial G}{\partial {\bf u}_c}(U),{\bf u}_b\right>+\left<\displaystyle\frac{\partial G}{\partial {\bf u}_b}(U),{\bf u}_c\right>\right),\end{array}\right.
		\end{equation*}
	where $1\leq a\leq p$, $1\leq b<c\leq p$. Form the symmetric matrix $\Sigma(U)$.
\item [3.] Choose the retraction $\mathcal{R}:T_{U}St^n_p\subset T_U\mathcal{M}_{n\times p}(\R)\rightarrow St^n_p\subset\mathcal{M}_{n\times p}(\R)$, $\mathcal{R}_U({\bf v}_U)=\text{qf}(U+{\bf v}_U)$, where $\text{qf}(A)$ denotes the $Q$ factor of the decomposition of $A \in \mathcal{M}_{n\times p}(\R)$
as $A=QR$, where $Q$ belongs to $St^n_p$ and $R$ is an upper triangular $n\times p$ matrix with strictly positive diagonal elements (see \cite{absil-mahony-sepulchre-1}).

\item [4.] Input ${U}^{(0)}\in St^n_p$ and $k=0$.

\item [5.] {\bf repeat}

$\bullet$  Determine a set $I_p(U^{(k)})$ containing the indexes of the rows that form a full rank submatrix of $U^{(k)}$.
Construct the basis $\mathcal{B}_{U^{(k)}}=\mathcal{B}_{U^{(k)}}^{'}\cup \mathcal{B}_{U^{(k)}}^{''}$ using the equations \eqref{delta-prim} and \eqref{delta-secund}.

$\bullet$ Compute the coordinate functions:
\begin{align*}
 g'_{ab}(U^{(k)})= & (-1)^{a+b}\left(\left<\displaystyle\frac{\partial G}{\partial {\bf u}_b}(U^{(k)}),{\bf u}_a^{(k)}\right>-\left<\displaystyle\frac{\partial G}{\partial {\bf u}_a}(U^{(k)}),{\bf u}_b^{(k)}\right>\right),\,\,\,1\leq a<b\leq p \\
 g''_{ic}(U^{(k)})= & \left<\displaystyle\frac{\partial G}{\partial {\bf u}_c}(U^{(k)}),{\bf z}_i^{(k)}\right>,\,\,\,i\in \{1,2,...,n\}\setminus I_p(U^{(k)}),\,c\in \{1,2,...,n\},
\end{align*}
where ${\bf z}^{(k)}_i$ is the vector formed with the $i$-th column of the matrix $Z^{(k)}=\mathbb{I}_n-U^{(k)}{U^{(k)}}^T$.

$\bullet$ Compute the components of the Hessian matrix $\text{Hess}\,\wt{G}$ of the cost function $\wt{G}$ (using formulas \eqref{sigma-I-0}, \eqref{sigma-I-1}, and \eqref{sigma-I-2}):
{\begin{align*}\label{Hessian-explicit}
 h_{\alpha_1\beta_1,\alpha_2\beta_2}(U^{(k)}) & = \left(\text{Hess}\, G(U^{(k)})-{\Sigma}(U^{(k)})\otimes \mathbb{I}_n\right)\left(\Delta'_{\alpha_1\beta_1}(U^{(k)}),\Delta'_{\alpha_2\beta_2}(U^{(k)})\right)\nonumber\\
  h_{\alpha\beta,jd}(U^{(k)}) & =\, \text{Hess}\, G(U^{(k)})\left(\Delta'_{\alpha\beta}(U^{(k)}),\Delta''_{jd}(U^{(k)})\right) \\
 h_{j_1d_1,j_2d_2}(U^{(k)}) & = \left(\text{Hess}\, G(U^{(k)})-{\Sigma}(U^{(k)})\otimes \mathbb{I}_n\right)\left(\Delta''_{j_1d_1}(U^{(k)}),\Delta''_{j_2d_2}(U^{(k)})\right),\nonumber
\end{align*}}
where $\alpha_1,\beta_1,\alpha_2,\beta_2,\alpha,\beta\in \{1,2,...,p\}$ with $\alpha_1<\beta_1$, $\alpha_2<\beta_2$, $\alpha<\beta$, and $j,j_1,j_2\in \{1,2,...,n\}\setminus I_p(U^{(k)})$ and $d,d_1,d_2 \in \{1,2,...,p\}$.

$\bullet$ Solve the linear system \eqref{Newton-equation-forms} with the unknowns  $\left(...,v^{\alpha\beta},...,v^{jd},...\right)$:
\begin{equation*}
\left[
\begin{array}{c|c}
h_{\alpha_1\beta_1,\alpha_2\beta_2}(U^{(k)})\, &\, h_{\alpha_1\beta_1,j_2d_2}(U^{(k)}) \\
\hline 
 h_{j_1d_1,\alpha_2\beta_2}(U^{(k)}) & h_{j_1d_1,j_2d_2}(U^{(k)}) 
\end{array}\right]
\cdot 
\left[
\begin{array}{c}
v^{\alpha_2\beta_2} \\
\hline
v^{j_2d_2}
\end{array}
\right]=
-
\left[
\begin{array}{c}
g'_{\alpha_1\beta_1}(U^{(k)}) \\
\hline 
g''_{j_1d_1}(U^{(k)})
\end{array}
\right].
\end{equation*}

\end{itemize}

$\bullet$ Construct the line search vector 
$$ {\bf v}_{U^{(k)}}=\sum\limits_{1\leq \alpha<\beta\leq p}  v^{\alpha\beta}\Delta'_{\alpha\beta}(U^{(k)})+\sum\limits_{
{\footnotesize \left.\begin{array}{l}
j\in \{1,2,...,n\}\setminus I_p(U^{(k)}) \\
d\in \{1,2,...,n\}
\end{array}
\right.}
} v^{jd}\Delta''_{jd}(U^{(k)}).$$

$\bullet$ Set $U^{({k+1)}}=\mathcal{R}_{U^{(k)}}\left({\bf v}_{U^{(k)}}\right)=\text{qf}\left(U^{(k)}+{\bf v}_{U^{(k)}}\right)$.

{\bf until} $U^{({k+1)}}$ sufficiently minimizes $\wt{G}$.

\end{framed}

For particular cases of cost functions defined on Stiefel manifolds, using ingenious matrix representations of the Hessian and the gradient vector field, in \cite{sato-1} and \cite{sato-2}  are given explicit formulas to solve the Newton equation.

\section{Conclusions}

In Theorem \ref{Hessian-matricial-123} we present a general formula for the Hessian matrix of a cost function defined on an orthogonal Stiefel manifold. We also point out the explicit expressions for the Lagrange multiplier functions, which are defined on the whole manifold and not just in the critical points of the cost function. This fact makes the formula for the Hessian matrix suitable for explicitly writing numerical algorithms like, for example, the Newton method, once we determine a basis for the tangent space in a given point of the orthogonal Stiefel manifold. Such a basis is explicitly constructed in Section \ref{section-3} and its properties are analyzed. We show that this construction decisively depends on the choice of a full rank submatrix for the given point in the Stiefel manifold.

\subsection*{Acknowledgements}
This work was supported by a grant of Ministery of Research and Innovation, CNCS - UEFISCDI, project number PN-III-P4-ID-PCE-2016-0165, within PNCDI III.\\
The authors would like to thank the anonymous referees for their valuable comments that helped improve the paper significantly.


\begin{thebibliography}{99}
\bibitem{absil-mahony-sepulchre-1} Absil, P.A., Mahony, R., Sepulchre, R.: Optimization Algorithms on Matrix Manifolds. Princeton University Press (2008)
\bibitem{sato-2}Aihara, K., Sato, H.: A matrix-free implementation of Riemannian Newton's method on the Stiefel manifold. Optim. Lett. 11, 1729-1741 (2017)
\bibitem{birtea-casu-comanescu} Birtea, P., Ca\c su, I., Com\u anescu, D.: First order optimality conditions and steepest descent algorithm on orthogonal Stiefel manifolds. Optim. Lett. 13, 1773-1791 (2019)
\bibitem{birtea-comanescu} Birtea, P., Com\u anescu, D.: Geometric dissipation for dynamical systems. Comm. Math. Phys. 316, 375-394 (2012)
\bibitem{Birtea-Comanescu-Hessian}Birtea, P., Com\u anescu, D.: Hessian operators on constraint manifolds. J. Nonlinear Science 25, 1285-1305 (2015)
\bibitem{5-electron} Birtea, P., Com\u anescu, D: Newton algorithm on constraint manifolds and the 5-Electron Thomson problem. J. Optim. Theor. Appl. 173, 563-583 (2017)
\bibitem{chu} Chu, M.T.,  Trendafilov, N.T.: The orthogonally constrained regression revisited.  J. Comput. and Graphical Statistics 10, 746-771 (2001)
\bibitem{dodig}Dodig, M., Sto\v{s}i\'{c}, M., Xavier, J.: On minimizing a quadratic function on Stiefel manifold, Linear Algebra and its Applications 475, 251-264 (2015)
\bibitem{edelman}Edelman, A., Arias, T.A., Smith, S.T.: The
geometry of algorithms with orthogonality constraints. SIAM J. Matrix Anal. Appl. 20, 303-353 (1998)
\bibitem{elden-park} Eld\'en, L., Park, H.: A Procrustes problem on the Stiefel manifold. Numer. Math. 82, 599-619 (1999)
\bibitem{gabay}Gabay, D.: Minimizing a differentiable function over a differential manifold. J. Optim. Theory Appl. 37, 177-219 (1982)
\bibitem{lu}Lu, J., Davidson, T.N.,  Luo, Z.-Q.: Blind separation of BPSK signals using Newton's method on the Stiefel manifold. IEEE International Conference on Acoustics, Speech, and Signal Processing 4, 301-304 (2003)
\bibitem{manton}Manton, J.H.: A framework for generalising the Newton method and other iterative methods from Euclidean space to manifolds. Numer. Math. 129, 91-125 (2015)
\bibitem{sato-1}Sato, H.: Riemannian Newton-type methods for joint diagonalization on the Stiefel manifold with application to independent component analysis. A Journal of Mathematical Programming and Operations Research 66, 2211-2231 (2017)
\bibitem{shukla}Shukla, A., Anand, S.: Metric learning based automatic segmentation of patterned species. IEEE Conference on Image Processing, Phoenix, AZ, USA, September 2016
\bibitem{turaga}Turaga, P., Veeraraghavan, A., Chellappa, R.: Statistical analysis on Stiefel and Grassmann manifolds with applications in computer vision. IEEE Conference on Computer Vision and Pattern Recognition, Anchorage, AK, USA, June 2008
\bibitem{wen}Wen, Z., Yin, W.: A feasible method for optimization with orthogonality constraints. Math. Progr. (Series A) 142, 397-434 (2013)


\end{thebibliography}
\end{document}